\documentclass[]{gENO2e}
\usepackage{algorithm}
\usepackage{algorithmic}
\usepackage{multirow}
\usepackage{booktabs}
\usepackage{rotating}
\usepackage{color}
\usepackage{threeparttable}

\begin{document}

\doi{10.1080/0305215X.YYYY.CATSid}
 \issn{1029-0273}
\issnp{0305-215X}
\jvol{00} \jnum{00} \jyear{2013} \jmonth{October}

\markboth{Peng Guo Wenming Cheng and Yi Wang}{Engineering Optimization}


\title{A modified  generalized extremal optimization algorithm for the quay crane scheduling problem with interference constraints}

\author{Peng Guo$^{a}$
\vspace{6pt} Wenming Cheng $^{a}$\vspace{6pt} and Yi Wang$^{b}$$^{\ast}$\thanks{$^\ast$Corresponding author. Email: ywang2@aum.edu}\\\vspace{6pt}  $^{a}${\em{School of Mechanical Engineering, Southwest Jiaotong University, Chengdu, China}};\\
$^{b}${\em{Department of Mathematics, Auburn University at Montgomery, Montgomery, AL, USA}}\\\vspace{6pt}}

\maketitle

\begin{abstract}
The quay crane scheduling problem (QCSP) determines the handling sequence of tasks at ship bays by a set of cranes assigned to a container vessel such that the vessel's service time is minimized.
   A number of heuristics or meta-heuristics have been proposed to obtain the near-optimal solutions to overcome the NP-hardness of the problem. In this article,
the idea of generalized extremal optimization (GEO) is adapted to solve the QCSP with respect to various interference constraints. The resulted GEO is termed as the modified GEO. A randomized searching method for neighboring task-to-QC assignments to an incumbent task-to-QC assignment is developed in executing the modified GEO. In addition, a unidirectional search decoding scheme is employed to transform a task-to-QC assignment to an active quay crane schedule.
The effectiveness of the developed GEO is tested on a suite of benchmark problems introduced by \citet{KimPark2004}. Compared with other well known existing approaches, the experiment results show that the proposed modified GEO is capable of obtaining the optimal or near-optimal solution in reasonable time, especially for large-sized problems.

\bigskip

\bigskip
\end{abstract}
\newpage

\vspace{-42pt}

\section{Introduction}

As the trend of globalization inevitably moves forward, the containerized maritime transportation has increased steadily in the past decades. According to the report of the United Nations Conference on Trade and Development, more than 80\% of the world freight is transported by sea. Port container terminals play a key role in serving as multi-modal interfaces between sea and land transportation. Thus container terminals are called to operate in the most efficient way, that is, as fast as possible at the least possible cost. The main operation of a container terminal is to process the loading and/or unloading tasks of vessels. These operations are performed by quay cranes (QCs). QCs are the most important equipments used at a terminal and their operational performance is a crucial factor of   turnaround times of vessels.

A well planned crane service schedule results in greatly reducing the vessel's service time (\emph{makespan}). A QC schedule specifies the service sequence of ship-bays of a container vessel by each QC and the time schedule for the services. The scheduling of cranes is typically constrained by many practical considerations. Since QCs along the same berth are sequentially mounted on two rail tracks, they cannot cross each other. This is referred to as \emph{non-crossing constraint}. For ensuring safe operations of cranes, \emph{safety margin} between adjacent cranes must be kept at any time. In general, in order to conveniently formulate the quay crane scheduling problem (QCSP), most researchers assumed that all containers in the same bay to be unloaded and/or loaded are treated as one single task. No preemption is allowed among all tasks. In addition, the earliest ready time of each crane and the traveling time for moving a crane from one bay to another bay are taken into consideration. Several studies have proposed methods to obtain the optimal schedule of QCs with consideration of above mentioned \emph{interference constraints}. However, most methods are not fit for solving a large-sized problem. Therefore, there is still a need of effective approaches to solve the QCSP.

The article focus on adapting the idea of generalized extremal optimization (GEO) developed by \citet{DeSousaRamos2003} to deal with the QCSP with respect to various interference constraints.
  The resulted GEO is termed as the \emph{modified GEO} (MGEO). A randomized searching method for neighboring task-to-QC assignments to an incumbent task-to-QC assignment is developed in executing the MGEO. In executing the MGEO, a decoding scheme based on the unidirectional movement of QCs is used to transform a task-to-QC assignment to a feasible schedule. As an outline for the remainder of the article, Section \ref{sec:Sec2} gives the literature review. Section \ref{sec:Sec3} describes the QCSP problem and formulation. Section \ref{sec:Sec4} presents
 the MGEO approach, followed by computational experiments in Section \ref{sec:sec5}. Finally, Section \ref{sec:sec6} gives some conclusion remarks and discusses some future work.

\section{Literature review\label{sec:Sec2}}

In the literature, the QCSP has received great attention over the decades. Two latest surveys about the berth allocation and the QCSP can be found in \citet{BierwirthMeisel2010,Rashidi20133601}. \citet{Daganzo1989} firstly presented the static and dynamic QCSP and suggested an algorithm for determining the number of cranes to be assigned to ship-bays of multiple vessels. Since then, the QCSP has drawn a worldwide attention as port container terminals are developing rapidly. \citet{PeterkofskyDaganzo1990} subsequently proposed a branch and bound algorithm for the static crane scheduling problem. \citet{LimRodrigues2004} augmented the study of Peterkofsky and Daganzo by incorporating  non-crossing constraints for a single vessel, and proposed some approximation algorithms. \citet{KimPark2004} defined a task as an unloading and loading operation for a cluster of adjacent slots on a deck in their study. They considered the non-crossing and precedence constraints for the QCs and formulated the problem as a MIP model. In their solution to the problem, a branch and bound method (B\&B) and a greedy randomized adaptive search procedure (GRASP) were designed for the solution.

The QCSP with complete bays introduced by \citet{Daganzo1989} and \citet{LimRodrigues2004} was studied by \citet{ZhuLim2006}, \citet{LimRodrigues2007}, \citet{LeeWang2008a}, \citet{LeeWang2008b} and \citet{LeeChen2010}. In \citet{ZhuLim2006}, the QCSP with non-crossing constraints (QCSPNC) was proven to be NP-complete and was solved by a branch and bound algorithm and a simulated annealing heuristic. In the work of \citet{LeeWang2008a}, a genetic algorithm (GA) was proposed to obtain near optimal solutions for the QCSPNC. Meanwhile, they used a GA to solve the QCSP with handling priority \citep{LeeWang2008b}. \citet{LeeChen2010} further proposed an approach to avoid the occurrence of unrealistic optimal solutions and designed two approximation schemes for the QCSPNC.

In solving the QCSP with container groups introduced by \citet{KimPark2004}, \citet{MocciaCordeau2006}  developed a branch and cut algorithm (B\&C) incorporating inequality constraints adopted from the solution methods for the precedence constrained traveling salesman problem. Subsequently, \citet{SammarraCordeau2007} proved that the QCSP can be viewed as a vehicle routing problem, which can be decomposed into a routing problem and a scheduling problem. They applied a tabu search (TS) to solve the routing problem and proposed a local search technique to generate the solution. Compared with the branch and cut algorithm and the GRASP, the TS algorithm provides a good balance between the solution quality and the computation time. Furthermore, \citet{NgMak2006} developed a scheduling heuristic to find effective schedules for the QCSP. \citet{BierwirthMeisel2009} proved that the unidirectional search (UDS) heuristic outperforms the standard solver CPLEX and some recent competing approaches \citep{KimPark2004,MocciaCordeau2006,SammarraCordeau2007}. \citet{Chung20124213} proposed a modified genetic algorithm to solve the QCSP. Their approach, compared to the results attained by other well known existing approaches, obtains better solutions in small sized instances and some medium sized instances. \citet{Meisel2011} introduced the QCSP with time windows, and presented a tree-search-based heuristic to solve the problem. Moreover, \citet{LegatoTrunfio2012} considered a rich QCSP that incorporates some practical considerations, and proposed a new mathematical formulation. In their work, a branch and bound scheme incorporating a unidirectional scheduling paradigm was employed.

The generalized extremal optimization (GEO) was originally developed as an improvement of the extremal optimization \citep{DeSousaRamos2003}. The GEO has been applied to some complex optimization problems \citep{DeSousaVlassov2004,DeSousaVlassov2004HTE,DeSousaSoeiro2007,GalskiDeSousa2007,Chen2007115,CassolSchneider2011,CucodeSousa2011}, such as the optimal design of a heat pipe, the spacecraft thermal design, the traveling salesman problem and so on. Recently, \citet{SwitalskiSeredynski2010,SwitalskiSeredynski2012} used the GEO to solve a multiprocessor scheduling problem. Their work indicated that the GEO demonstrates better average performance compared to the  GA and the particle swarm optimization. \emph{To the best knowledge of the authors, there is no literature that uses the GEO to solve the QCSP. Therefore, it is intended in this article to adapt a GEO approach to obtain near-optimal solutions for the QCSP, while interference constraints are considered.}

\section{Problem description and formulation\label{sec:Sec3}}

\begin{figure}[h]
  \centering
  \includegraphics[scale=0.7]{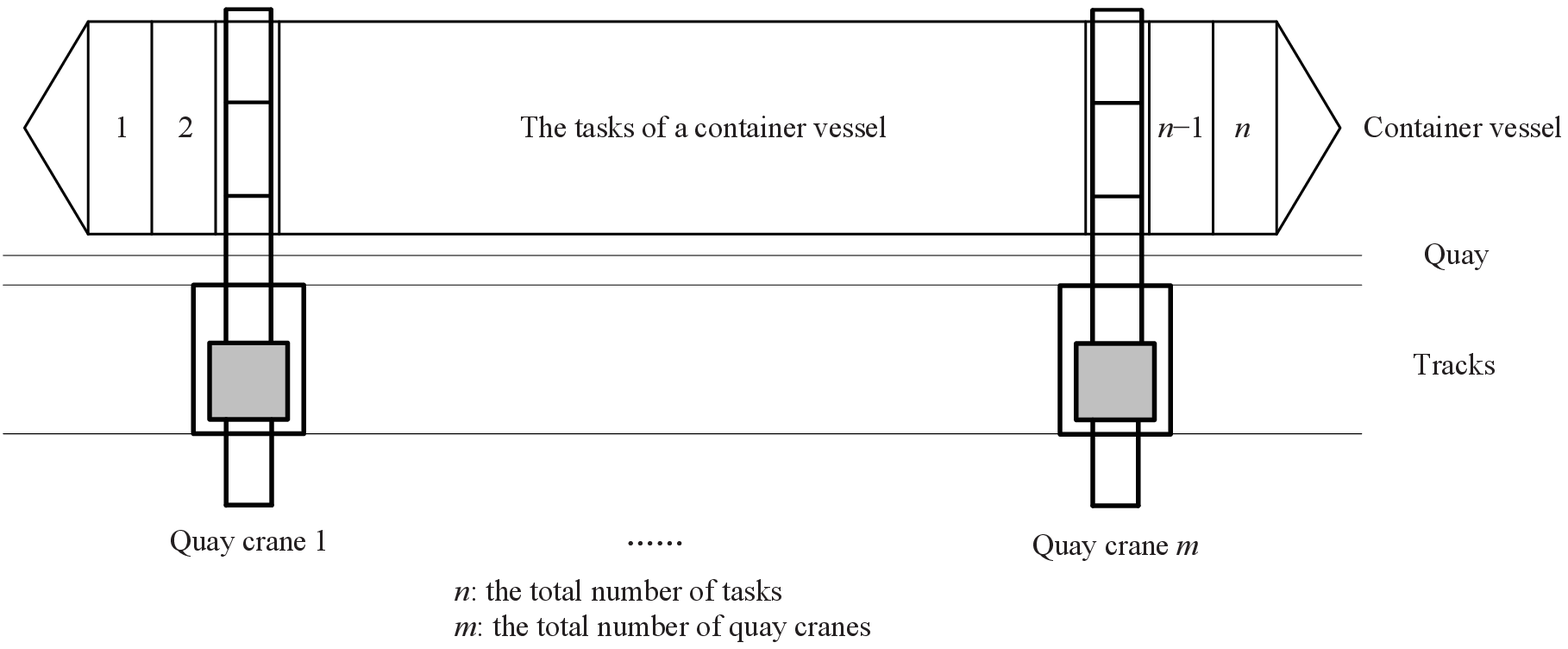}\\
  \caption{An illustration of QCs working on a container vessel with $l=n$.}\label{fig: fig_1}
\end{figure}

The considered QCSP focuses on a single freight vessel with $l$ bays where $n$ tasks must be processed by a set of identical quay cranes $Q=\{1, 2, \ldots, m\}$ that are allocated to the vessel. Let $\Omega$ denote the set of the $n$ tasks. Each task $i\in \Omega$ represents a loading/unloading operation of a certain container group. Task $i$ has an individual processing times $p_i$ and its bay position is indicated by $l_i$. The bays of the vessel are labeled sequentially from bay 1 to bay $l$. The cranes are subsequently indexed in an ascending order along the direction of increasing bay positions, as shown in figure \ref{fig: fig_1}. Each crane $k\in Q$ has an earliest ready time $r_{0}^{k}$ and an initial bay position $l_{0}^{k}$. Since no two QCs can process different tasks at a bay simultaneously, pairs of tasks located in the same bay have \emph{precedence relations}. Let $\Phi$ represent the set of task pairs with precedence relations. For each task pair $(i, j)\in \Phi$, task $i$ is completed before task $j$ starts. The set $\Psi$, which is obtained in advance, contains pairs of tasks that require \emph{non-simultaneous processing}. Specifically, for each task pair $(i, j)\in \Psi$, either task $i$ must be completed before task $j$ starts or task $i$ starts only after task $j$ is completed. Clearly, $\Psi \supseteq \Phi$.

All cranes are mounted on the same tracks. They can move from one bay to an adjacent bay within a uniform traveling time $t_0$. The traveling time $t_{0i}^{k}$ of crane $k$ from its initial location to the bay position of task $i$ is calculated as $t_{0i}^{k}=|l_{i}-l_{0}^{k}|\cdot t_{0}$. However, keep in mind that QCs can not cross each other, and to ensure safety during the handling process, adjacent QCs must keep a safety margin at any time. In order to avoid violating the non-crossing and the safety margin constraints, a \emph{minimum travel time} $\Delta_{ij}^{vw}$ formulated in \citet{BierwirthMeisel2009} and \citet{Meisel2011} is inserted in a schedule between the processing of tasks $i$ and $j$ that are assigned to cranes $v$ and $w$, respectively, should interference of cranes occur.

Specifically, let $\delta$ be the minimal safety margin between adjacent cranes, usually expressed by a number of bays. It is easy to see that the smallest allowed difference $\delta_{vw}$ between the bay positions of QCs $v$ and $w$ is given by $\delta_{vw}=(\delta+1)\cdot |v-w|$. Consequently, the quantity $\Delta_{ij}^{vw}$ can be calculated as follows.
\begin{equation}\label{eq:eq1}
  \Delta_{ij}^{vw}=
  \begin{cases}
  (l_{j}-l_{i}+\delta_{vw})\cdot t_{0}, & \text{if $v > w$ and $i \neq j$ and $l_{i} < l_{j}+\delta _{vw}$,} \\
  (l_{i}-l_{j}+\delta _{vw})\cdot t_{0}, & \text{if $v < w$ and $i \neq j$ and $l_{i} > l_{j}-\delta_{vw}$,} \\
  |l_{i}-l_{j}|\cdot t_{0}, & \text{if $v=w$ and $i \neq j$,} \\
  0, & \text{otherwise.}
  \end{cases}
\end{equation}
In equation \eqref{eq:eq1}, the first case considers the situation in which crane $v$ operates above crane $w$, that is, $v>w$. If task $i$ that is processed by crane $v$ lies within the safety margin of task $j$, the two tasks cannot be processed simultaneously. Consequently, a minimum travel time $(l_{j}-l_{i}+\delta_{vw})\cdot t_{0}$ must be inserted between their processing to avoid the conflict of cranes. The reverse positioning situation of cranes $v$ and $w$ is tackled in the second case of equation \eqref{eq:eq1}. Besides, the situation where both tasks $i$ and $j$ are processed by the same crane $v=w$ is dealt with in the third case of equation \eqref{eq:eq1}. In this case the value of $\Delta_{ij}^{vw}$ is just the traveling time for the crane to move from bay position $l_{i}$ to bay location $l_{j}$. In all other cases, cranes do not get into a conflict, and there is no need to insert extra travel time as described in equation \eqref{eq:eq1} by setting  $\Delta_{ij}^{vw}=0$. For convenience the set $\Theta=\{(i, j, v, w) \in \Omega^{2} \times Q^{2} | i<j, \Delta_{ij}^{vw}>0\}$ that contains all combinations of a pair of tasks and a pair of cranes that potentially lead to crane interference is defined. \emph{Note that} if $(i, j, v, w)\in \Theta$,  {the task pair} $(i, j) $\emph{is not      subject to the precedence constraint}. \emph{Inserting the quantity $ \Delta_{ij}^{vw}$ effectively prohibits that task} $i$ and $j$ \emph{are serviced simultaneously}. Thus if task $j$ starts before task $i$ is completed, then the latter task $j$ should start at $s_{j} \geqslant c_{i}+\Delta_{ij}^{vw}$; or if task $i$ is processed before task $j$ is completed, then the latter task $i$ should start at $s_{i} \geqslant c_{j}+\Delta_{ij}^{vw}$. Consequently, \emph{an insertion of a positive $\Delta_{ij}^{vw}$ automatically handles the constraints given in sets $\Phi$ and $\Psi$}.

The mathematical model used for the considered QCSP is
based on the one proposed by \citet{KimPark2004} and the formulation of $\Delta_{ij}^{vw}$ in \citet{BierwirthMeisel2009} and \citet{Meisel2011}. The model's objective is to minimize the vessel's service time $C_{\max}$. In the model, the binary variable $x_{i}^{k}$ equals to 1 if and only if task $i \in \Omega$ is processed by crane $k \in Q$, and to 0 otherwise; the variable $y_{ij}$ is equal to 1 if and only if task $i\in \Omega$ is completed no later than task $j\in \Omega$ starts, and to 0 otherwise; the variable $c_{i}$ represents the completion time of task $i\in \Omega$; and the constant $M$ is some large positive number such that $M >\displaystyle \sum_{i \in \Omega} p_{i}$. The model is described in equations  \eqref{eq:3.2}--\eqref{eq:3.13}.
\enskip

minimize
\begin{eqnarray}
 C_{\max}  \label{eq:3.2}
\end{eqnarray}
subject to
\begin{align}
c_{i} & \leqslant C_{\max} & \forall i & \in \Omega \label{eq:3.3} \\
\sum_{k \in Q} x_{i}^{k} & =1 & \forall i & \in \Omega \label{eq:3.4} \\
r_{0}^{k}+t_{0i}^{k}+p_{i}-c_{i} & \leqslant M \cdot (1-x_{i}^{k}) & \forall i & \in \Omega, k \in Q \label{eq:3.5} \\
c_{i}-c_{j}+p_{j} & \leqslant 0 & \forall (i,j) & \in \Phi \label{eq:3.6} \\
y_{ij}+y_{ji} & =1 & \forall (i,j) & \in \Psi \label{eq:3.7} \\
c_{i}+p_{j}-c_{j} & \leqslant M \cdot (1-y_{ij}) & \forall i,j &\in \Omega \label{eq:3.8} \\
c_{j}-p_{j}-c_{i} & \leqslant M \cdot y_{ij} & \forall i,j &\in \Omega \label{eq:3.9} \\
x_{i}^{v}+x_{j}^{w} & \leqslant 1+y_{ij}+y_{ji} & \forall (i,j,v,w) & \in \Theta \label{eq:3.10} \\
c_{i}+\Delta_{ij}^{vw}-c_{j}+p_{j} & \leqslant M \cdot (3-y_{ij}-x_{i}^{v}-x_{j}^{w}) & \forall (i,j,v,w) & \in \Theta \label{eq:3.11} \\
c_{j}+\Delta_{ij}^{vw}-c_{i}+p_{i} & \leqslant M \cdot (3-y_{ji}-x_{i}^{v}-x_{j}^{w}) & \forall (i,j,v,w) & \in \Theta \label{eq:3.12} \\
x_{i}^{k},y_{ij} & \in \{0,1\} & \forall i,j & \in \Omega, k \in Q \label{eq:3.13}
\end{align}

In the above formulation, constraints \eqref{eq:3.3} define the property of $C_{\max}$. Constraints \eqref{eq:3.4} ensure that every task must be assigned to exactly one crane. Constraints \eqref{eq:3.5} give the correct completion time for the first task of each crane. Precedence relations between tasks in the set $\Phi$ are respected through Constraints \eqref{eq:3.6}. Constraints \eqref{eq:3.7} ensure that tasks $i$ and $j$ cannot be processed simultaneously if $(i, j) \in \Psi$. Constraints \eqref{eq:3.8} and \eqref{eq:3.9} determine the completion times of tasks with respect to the variables $y_{ij}$. The conflicts among cranes are considered in constraints \eqref{eq:3.10}-\eqref{eq:3.12}. Constraints \eqref{eq:3.10} guarantee that each pair of tasks $i$ and $j$ from the set $\Theta$ are not processed simultaneously if they are assigned to cranes $v$ and $w$, respectively. In constraints \eqref{eq:3.10}, $x_{i}^{v}=1$ means that task $i$ is assigned to crane $v$ and $x_{j}^{w}=1$ means that task $j$ is assigned to crane $w$. If these assignments come together, tasks $i$ and $j$ cannot be processed simultaneously, which implies that $y_{ij}$ and $y_{ji}$ are not set to 1 at the same time. The case of $y_{ij}=1$ is considered in constraints \eqref{eq:3.11}. The case $y_{ji}=1$ is tackled in constraints \eqref{eq:3.12}. The two constraints insert the necessary minimum travel time $\Delta_{ij}^{vw}$ between the processing of tasks $i$ and $j$ for cases $y_{ij}=1$ and $y_{ji}=1$, respectively. Finally, constraints \eqref{eq:3.13} define the ranges for the decision variables.

\section{The modified generalized extremal optimization algorithm\label{sec:Sec4}}

The QCSP can be viewed as a scheduling problem on parallel identical machines with various constraints. It has been proven that the QCSP is NP-hard \citep{SammarraCordeau2007}. Therefore, there exists no polynomial time algorithm to exactly solve the considered problem. In this section, a modified generalized extremal optimization approach based on the unidirectional schedule is designed to obtain near-optimal solutions. During the search process of the proposed algorithm, both non-crossing constraints and safety margin constrains are considered.
The generalized extremal optimization (GEO) is a recently proposed meta-heuristic devised to solve complex optimization problems \citep{DeSousaRamos2003}. It was originally developed as a generalization of the extremal optimization approach \citep{BoettcherPercus2001}, which was inspired by the evolution model of Bak-Sneppen \citep{BakSneppen1993}. It is of easy implementation with only one free parameter to be adjusted, and does not make use of derivatives. The algorithm can be applied to a broad class of nonlinear constrained optimization problems.

In the GEO algorithm, a population of species is represented by a binary string of $L$ bits that encodes the design variables of the optimization problem.   Each bit of a string indicates a species. Each bit is associated to a fitness number that is proportional to the gain or loss of the cost function value if the bit is mutated (changed from 1 to 0 or vice versa). Then, all bits are ranked by their fitness numbers from $\kappa=1$, for the least adapted bit, to $\kappa=L$ for the best adapted. A bit is selected to mutate according to the probability distribution $P_{\kappa}\propto \kappa^{-\tau}$, where $\kappa$ is the rank of a randomly chosen bit candidate to mutate and $\tau$ is a positive parameter. If $\tau \rightarrow 0$, all bits have the same probability to mutate, whereas if $\tau \rightarrow \infty$ the worst adapted bit will mutate.

It has been proven that the GEO is suitable for a task scheduling problem due to its simplicity \citep{SwitalskiSeredynski2010}. In what follows, \emph{the idea of the GEO is applied to develop an algorithm termed as the modified GEO (MGEO) to solve the QCSP, meanwhile considering the non-crossing and safety margin constrains during the search process}.

\emph{In the MGEO, implemented is a unidirectional schedule (UDS) for the order of cranes visiting the bays}. In a UDS, the QCs do not change the moving direction after possible initial repositioning and have identical moving directions either from lower to upper bays or vice versa. Generally, to find the best UDS one needs to consider the case in which the cranes move from lower to higher bay positions but also the reverse case. Therefore, the decoding procedure needs to apply twice to find the best UDS. First, a UDS is generated for an \emph{upward movement} of QCs. To generate a UDS in the reverse direction, the bays are instead numbered in reverse order. Accordingly, the tasks and QCs are renumbered in the order of increasing bays. The precedence constraints, non-simultaneous tasks are updated with respect to the new task numbering. The initial positions of QCs are subsequently adapted as well. Then the same UDS with upward movement of QCs is applied to generate a best schedule, which is actually the best unidirectional schedule for a \emph{downward moment} of QCs. Consequently, the UDS delivers the better of the two schedules.

Therefore, in what follows, an \emph{upward movement} of QCs is fixed to build a schedule for a given task-to-QC assignment. The notation $a=(q_{1},q_{2},\ldots , q_{n})$ is introduced to represent a particular task-to-QC assignment sequence, where $q_{j},\ j \in \Omega$ is the crane assigned to task $j$. The initial task-to-QC assignment sequence shall be denoted by $a_{0}$. Furthermore, assume that the task-to-QC assignment to crane $k, k \in Q$ is denoted by the sequence $a^{k}=(k_{1},k_{2}, \ldots , k_{n_{k}})$,  where, the relation $k_{1}<k_{2}<\cdots <k_{n_{k}}$ must be true due to the adopted ordering of tasks and the assumption of an upward movement of cranes. For convenience, the notation $k_{0}$ is used to denote a \emph{dummy} task indicating the beginning of service of crane $k$. Observe that there are $m$ such sequences, $a^{1}$, $a^{2}$, $\ldots$, $a^{m}$.

\subsection{String representation of a task-to-QC assignment}

\begin{figure}[h]
  \centering
  \includegraphics[scale=1]{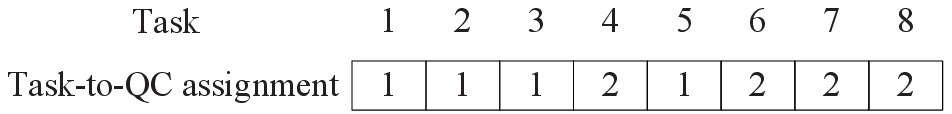}\\
  \caption{An illustration of a string representation for a task-to-QC assignment sequence.}\label{Fig:fig2}
\end{figure}

An admissible task-to-QC assignment sequence is represented in the MGEO by a string whose length is equal to the total number of tasks. The $i$-th bit of the string contains the crane to which the $i$th task is assigned. Figure \ref{Fig:fig2} shows a string representation for a task-to-QC assignment with eight tasks and two QCs.

\subsection{Initial task-to-QC heuristic}

To generate an initially feasible task-to-QC assignment for a QCSP instance, two approaches introduced by \citet{SammarraCordeau2007} have been utilized: the \emph{S-TASKS} and the \emph{S-LOAD}. The \emph{S-TASKS} approach distributes tasks to cranes so that each crane is assigned approximately equal number of tasks; while the \emph{S-LOAD} approach assigns tasks to all cranes so that the total processing time of tasks assigned to each crane is roughly even. According to the preliminary experiments, the best initial solution is obtained by the \emph{S-LOAD} approach. The procedure of  the \emph{S-LOAD} approach is described in figure \ref{alg:alg_SLOAD}. Thus, for the given instance in table \ref{tab:tab1}, the initial task-to-QC assignment is obtained using the \emph{S-LOAD} rule such that tasks 1-5 are assigned to crane 1 and tasks 6-8 are assigned to crane 2. The corresponding schedule can be obtained through the decoding procedure described in Subsection \ref{sec:decoding}.   The resulted makespan of the container vessel is 666 time units.

\begin{figure}[h]
\begin{algorithmic}[1]
\STATE input the initial data of a given instance;
\STATE consider the task set $\Omega$;
\STATE $p\_m=\sum_{i \in \Omega}p_{i}/m$;   \enskip \COMMENT{calculate the average of the processing times of all tasks}
\STATE set $a^{k}=[\ ]$ for $k\in Q$;   \enskip  \COMMENT{assume that the task-to-QC assignment to crane $k\in Q$ is empty}
\STATE $t^{k}=0$ for $k\in Q$;   \enskip  \COMMENT{initialize the completion times of all cranes}
\STATE $k=1$;   \enskip  \COMMENT{prepare to assign tasks to crane 1}
\FOR { $i\leftarrow1$ to $n$ }
\STATE $t^{k}=t^{k}+p_{i}$;
\IF {$t^{k}\leqslant p\_m$ and $k\leqslant m$}
\STATE $a^{k}=[a^{k} \  i]$;    \enskip \COMMENT{append task $i$ to the current task list $a^{k}$ for crane $k$}
\ELSIF {$t^{k}>p\_m$ and $k\leqslant m$}
\STATE $a^{k}=[a^{k} \  i]$;   \enskip \COMMENT{append the last task $i$ to crane $k$}
\STATE $k=k+1$;    \enskip \COMMENT{prepare to assign tasks to next crane}
\ELSIF {$k\geqslant m$}
\STATE break   \enskip  \COMMENT{complete the task-to-QC assignments for all cranes and terminates the execution of the for loop}
\ENDIF
\ENDFOR
\STATE output the task-to-QC assignment $a^{k}$ for $k \in Q$.
\end{algorithmic}
\caption{\label{alg:alg_SLOAD}The pseudo-code of the \emph{S-LOAD} approach.}
\end{figure}

\begin{table}[h]\scriptsize
  \centering
  \caption{A QCSP instance.}
    \begin{tabular}{rrrrrrrrr}
    \toprule
   \multicolumn{1}{l}{ Task ($i$)} & \multicolumn{1}{c}{1} & \multicolumn{1}{c}{2} & \multicolumn{1}{c}{3} & \multicolumn{1}{c}{4} & \multicolumn{1}{c}{5} & \multicolumn{1}{c}{6} & \multicolumn{1}{c}{7} & \multicolumn{1}{c}{8} \\
    \midrule
    \multicolumn{1}{l}{Processing time $p_{i}$} & \multicolumn{1}{c}{55} & \multicolumn{1}{c}{121} & \multicolumn{1}{c}{70} & \multicolumn{1}{c}{129} & \multicolumn{1}{c}{134} & \multicolumn{1}{c}{143} & \multicolumn{1}{c}{98} & \multicolumn{1}{c}{43} \\
   \multicolumn{1}{l}{ Bay position $l_{i}$} & \multicolumn{1}{c}{1} & \multicolumn{1}{c}{1} & \multicolumn{1}{c}{2} & \multicolumn{1}{c}{4} & \multicolumn{1}{c}{5} & \multicolumn{1}{c}{5} & \multicolumn{1}{c}{7} & \multicolumn{1}{c}{8} \\
  \multicolumn{1}{l}{QC 1} & \multicolumn{8}{l}{$l_{0}^{1}=$4, $r_{0}^{1}=$0} \\

   \multicolumn{1}{l}{QC 2} & \multicolumn{8}{l}{$l_{0}^{2}=$8, $r_{0}^{2}=$0} \\

   \multicolumn{1}{l}{ Precedence relations $\Phi$} & \multicolumn{8}{l}{\{(1, 2), (5, 6)\}} \\
  \multicolumn{1}{l}{  Non-simultaneous processing tasks $\Psi$} & \multicolumn{8}{l}{\{(1, 2), (1, 3), (2, 3), (4, 5), (4, 6), (5, 6), (7, 8)\}} \\
    \multicolumn{9}{l}{Number of tasks $n$=8; number of cranes $m$=2; } \\
    \multicolumn{9}{l}{safety margin $\delta$=1(bay); crane traveling speed $t_{0}=1$(time unit/bay).} \\
    \bottomrule
    \end{tabular}%
  \label{tab:tab1}%
\end{table}%

\subsection{Decoding scheme} \label{sec:decoding}

Once a task-to-QC assignment is obtained, a decoding scheme is applied to build a schedule for each task. Let $\Omega_{k}$ ($k \in Q$) be the set of tasks that are assigned to crane $k$ ($k \in Q$). Let the size of $\Omega_{k}$ be denoted by $n_k$. Thus $\displaystyle \sum_{k=1}^{m}n_{k}=n$. The decoding scheme mainly determines the correct starting time $s_{k_{j}}$ for every task $k_{j}\in \Omega_{k}$ that is assigned to QC $k$. The corresponding completion time for task $k_j$ is accordingly calculated by the equation $c_{k_{j}}=s_{k_{j}}+p_{k_{j}}$, $1\leqslant j \leqslant n_k$. The starting time $s_{k_{j}}$ of task $k_j \in \Omega_{k}$ ($1\leqslant j\leqslant n_k$) is simply given by the equation
\begin{equation}\label{eq:eq14}
  s_{k_{j}}=\max\left\{ c_{k_{j-1}}+\Delta_{k_{j-1},k_{j}}^{k,k}, \max_{\substack{q_{i}\in \Omega_{k+1} \\ (q_{i},k_{j},k+1,k) \in \Theta}}\left\{ c_{q_{i}}+\Delta_{q_{i},k_{j}}^{k+1,k} \right\} \right\}, \quad k\in Q, 1\leqslant j \leqslant n_k
\end{equation}
When $k=m$, the second part on the right-hand side is ignored. The first part in the formula computes a starting time for task $k_j$  to be the sum of the completion time of its immediately preceding task $k_{j-1}$ with a traveling time of crane $k$ from bay $l_{k_{j-1}}$ to bay $l_{k_{j}}$. In the case $j-1=0$, that is, $j=1$ and  $k_1$ is the first task assigned to crane $k$, define $\Delta_{k_{0},k_{1}}^{k,k}:=t_{0,k_1}^{k}=|l_0^k-l_{k_1}|\cdot t_0$ and the completion time $c_{k_0}=r_0^k$, that is, the earliest ready time of crane $k$.

The second part adjusts the first part when there is a crane $k+1$ operates above crane $k$, and there is an interference of cranes $k$ and $k+1$ had they performed their assigned jobs simultaneously, that is, $\Delta_{q_i,k_j}^{k+1,k}>0$. In this case, \emph{priority must be given to  crane $k+1$ with higher bay position due to the upward movement of QCs}. Thus, whenever $\Delta_{q_i,k_j}^{k+1,k}>0$, crane $k$ first has to wait until the upper crane $k+1$ completes its job, then a minimum travel time $\Delta_{q_i,k_j}^{k+1,k}$ is inserted to allow the upper crane to move away upwardly so that the interference of cranes $k$ and $k+1$ is avoided. Giving priority to an upper crane is because otherwise crane $k$ must change its moving direction to allow the upper crane $k+1$ to process  task $q_i$, which, of course, contradicts to the unidirectional search prerequisite. The priority given to an upper crane should interference between cranes occurs also dictates that \emph{the proposed algorithm necessarily starts with scheduling tasks for crane m at the uppermost bay position, followed by the next crane sequentially till crane 1. Moreover, the scheduling of tasks assigned to crane m is nearly determined by the first part of  equation} \eqref{eq:eq14}.

\subsection{Generating neighboring task-to-QC assignments}\label{subsec44}

Suppose that an incumbent task-to-QC assignment $a_c$ is configured and its corresponding unidirectional schedule is obtained according to the decoding scheme, to further execute the MGEO algorithm for a better schedule,  `\emph{neighboring task-to-QC assignments}' relative to $a_{c}$ are generated. To this end, let $a_c=(q_1, q_2, \ldots, q_n)$. For a task $j,\  1\leqslant j\leqslant n$, let the set $Q_j$ contains all cranes that can possibly service task $j$. Task $j$ (or bit $j$) is said to be mutated if the incumbent crane $q_j$ is replaced by a distinct crane $q_j^* \in Q_j$. Two or more tasks (to which shall be referred as multiple tasks hereafter) may be chosen to be mutated simultaneously in order to have bigger searching space for a better schedule. To describe this idea, denote by $\Omega_*$ the set of tasks that shall be simultaneously mutated for which neighboring task-to-QC assignments of $a_c$ shall be searched. Let $n_*$ be the size of $\Omega_*$. If $n_*=1$ then each task $j\in \Omega$ shall be assigned to $\Omega_*$ sequentially one at a time, that is, each time one task is mutated followed by the search of all its neighboring task-to-QC assignments. If $n_*\geqslant 2$, multiple tasks are chosen to be assigned to the set $\Omega_*$ simultaneously. Specifically, \emph{each time the multiple tasks to be mutated are chosen in a random way such that their indices correspond to distinct pseudo-random integers on the interval} [1, $n$]. \emph{Randomly choosing multiple tasks to be mutated avoids the searching for a better schedule trapped in a neighborhood of single locally optimized solution}.

After each task in $\Omega_*$ is mutated, each resulted new task-to-QC assignment is called a neighboring task-to-QC assignment (or a mutation) to $a_c$, denoted by $a_*$. The set of all the possible neighboring task-to-QC assignments relative to a specific $\Omega_*$ shall be generated by substituting every distinct crane $q_{j}^{*}\in Q_j$ for the incumbent crane $q_j$ for each task $j\in \Omega_*$. The resulted set of all neighboring task-to-QC assignments relative to $\Omega_*$ is denoted by $A_*$. If $\Omega_*$ contains only one task, that is, $n_*=1$, then the algorithm is run $n$ times so that in the $j$-th run, $1\leqslant j \leqslant n$, task $j\in \Omega$ is mutated to generate all of its neighboring task-to-QC assignments. Otherwise, if $\Omega_*$ contains two or more tasks, that is, $n_*\geqslant 2$, then the random procedure of choosing multiple tasks will  be repeated only $n$ times, where in each instance of $\Omega_*$, the set $A_*$ of all its neighboring task-to-QC assignments are generated.


Some practical considerations must be taken into consideration in order to guarantee a realistic mutation. As described in \citet{LeeChen2010}, due to the safety margin constrains, the space between two QCs may not be enough to allow one of them to operate a task positioned between them, or a QC may be incorrectly scheduled to be driven out of the boundary of the rail tracks during the effort to maintain the safety margin. {\em The key is to decide the set $Q_j$ for task $j$}. For a given task $j$ located in bay position $l_j$, the following rule is proposed to decide the set $Q_j$ of cranes that can possibly service task $j$.

Assuming the lowest bay position is 1, the highest bay position is $l$. If crane $q$, $q\in Q$, can service task $j$, due to the safety margin requirements between cranes, there must be sufficient space to hold the cranes in order to avoid incorrectly scheduling cranes to be driven out of the boundary of the rail tracks. Thus there must hold true that
\begin{equation}\label{eq:eq15}
  \delta_{1q}\leqslant l_{j}-1
\end{equation}
and
\begin{equation}\label{eq:eq16}
  \delta_{qm}\leqslant l-l_{j} \text{.}
\end{equation}

Since $1\leqslant q\leqslant m$, in view of formulas \eqref{eq:eq15} and \eqref{eq:eq16} the following two formulas are deduced:
\begin{equation}\label{eq:eq17}
  1\leqslant q+(q-1)\delta \leqslant l_{j} \leqslant m+(m-1)\delta
\end{equation}
and
\begin{equation}\label{eq:eq18}
  l-(m-1)(1+\delta)\leqslant l_{j} \leqslant l-(m-q)(1+\delta)\leqslant l \text{.}
\end{equation}
From equation \eqref{eq:eq17}, one can deduce that when $l_j\leqslant m+(m-1)\cdot\delta$, crane $q\in Q_j$ if $q$ is an integer and $1\leqslant q\leqslant (l_j+\delta)/(1+\delta)$; while by equation \eqref{eq:eq18} when $l_j\geqslant l-(m-1)(1+\delta)$, crane $q\in Q_j$ if $q$ is an integer and $m-(l-l_j)/(1+\delta)\leqslant q\leqslant m$. For all other cases, the set $Q_j$ is given by setting $Q_j=Q$.

\begin{figure}[h]
  \centering
  \includegraphics[scale=1]{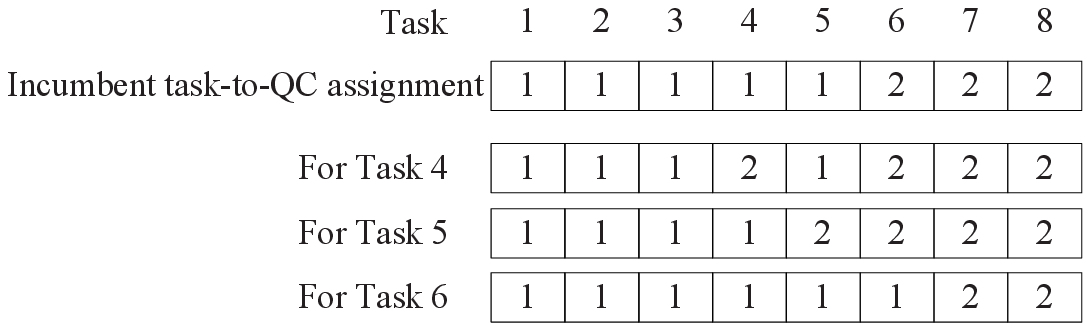}\\
  \caption{Example of neighboring task-to-QC assignments.}\label{fig:fig5}
\end{figure}

In figure \ref{fig:fig5}, all the neighboring task-to-QC assignments of a current incumbent task-to-QC assignment is displayed. The incumbent is the initial task-to-QC assignment obtained by the \emph{S-LOAD} rule for the instance in table 1. The neighboring task-to-QC assignments for tasks 1, 2, 3, 7, 8 are not displayed, as tasks 1, 2, 3 are only assigned to crane 1, and tasks 7, 8 are only assigned to crane 2, owing to the constraints of their bay positions governed by equations \eqref{eq:eq17} and \eqref{eq:eq18}.

\subsection{The proposed algorithm for the QCSP}

The proposed GEO algorithm starts with initiating a solution as the incumbent task-to-QC assignment and then launches a neighborhood search procedure. The procedure generates the set $A_*$ of all the neighboring task-to-QC assignments  based on the incumbent  $a_{c}$ by the method described  in subsection \ref{subsec44}. Let $n_{A_{*}}$ be the number of assignments in $A_*$. The $n_{A_{*}}$ assignments are ranked by their objective values, from $\kappa=1$ for the least one to $\kappa=n_{A_{*}}$ for the largest one.
A uniformly distributed random number $RAN$ on the interval [0,1] is generated. Then one neighboring assignment $a_{\kappa_i}(\kappa_i=1,\ldots,n_{A_{*}})$ is randomly chosen.  Let $\tau$ be a suitable positive constant that shall be chosen by experiments in terms of both solution quality and computing time. If the {\em updating probability} $P(\kappa_i)=\kappa_{i}^{-\tau}$ of the chosen assignment is equal to or greater than $RAN$, the neighboring assignment $a_{\kappa_i} $ is confirmed to be  chosen as the incumbent one. Otherwise, the process is repeated until a random  assignment in the set $A_*$ is chosen to be the incumbent one. Note that   in the proposed GEO algorithm the process of updating the incumbent task-to-QC assignment $a_c$ for the purpose of generating the set $A_*$ of  neighboring task-to-QC assignments allows the algorithm to search broader solution space. Consequently the proposed algorithm can  avoid the search process to get trapped in a local optimum and hence may obtain a better result.

Finally, the neighboring assignment $a_{\kappa_1}$ is compared to the current best assignment in terms of the objective value. If $a_{\kappa_1}$ is a better task-to-QC assignment, $a_{\mathrm{best}}$ is replaced by $a_{\kappa_1}$.
The search process iterates in a new neighborhood of the updated incumbent $a_c$   until the {\em stopping criterion} is met.
In this article, the stopping criterion when testing all instances is set to be the maximum iterations $I_{\max}$ fixed to 200. Meanwhile,   another stopping rule is used to reduce the computational time, that is,  the maximum number $I'_{\mathrm{max}}$ of  iterations with no-improvement is set to 50.
A pseudo-code of the proposed algorithm is shown in figure \ref{alg:GEO}.

\begin{figure}[ht]
 \begin{algorithmic}[1]
 \STATE \textbf{Initialization}:
 \STATE \enskip Obtain an initial task-to-QC assignment sequence $a_0$;
  \STATE \enskip Calculate the objective value $f(a_0)$ of the initial task-to-QC assignment $a_0$. Set the best solution $a_{\mathrm{best}}=a_0$ and the best solution value $f(a_{\mathrm{best}})=f(a_0)$;
  \STATE \enskip Choose the stopping criterion: $I_{\max}=200$, $I'_{\mathrm{max}}=50$;
  \STATE \enskip Set the incumbent task-to-QC assignment $a_c=a_0$. Initialize the counters $I=0$ and $I'=0$; Set $\tau=5$;
  \STATE For each $j$, $j\in \Omega$, determine the set $Q_j$ of cranes that can provide service to task $j$;
  \WHILE {the stopping criterion is not met, that is, when $I\leq I_{\max}$ or $I'\leq I'_{\mathrm{max}}$}

  \STATE For the incumbent $a_c$, choose $n_*$ and determine the set $\Omega_*$ of tasks that shall be simultaneously mutated;

  \STATE Generate the set $A_*$ of all possible neighboring task-to-QC assignments relative to $\Omega_*$;
  \STATE Calculate the objective function values of these assignments in the set $A_*$;
  \STATE Rank the $n_{A_{*}}$ assignments according to their objective function values;

  \REPEAT

  \STATE Generate an uniformly distributed random number $RAN\in [0,1]$;
  \STATE Choose randomly one neighboring assignment $a_{\kappa_i}$ from the set $A_*$;
  \IF {$RAN\leq P(\kappa_i)=\kappa_{i}^{-\tau}$}
  \STATE $a_c=a_{\kappa_i}$;
  \ENDIF

  \UNTIL {the incumbent $a_c$ is updated by a neighboring one $a_{\kappa_i}$ in the set $A_*$}
  \IF {$f(a_{\kappa_1})<f(a_{\mathrm{best}})$}
  \STATE $f(a_{\mathrm{best}})=f_{a_{\kappa_1}}$;
  \STATE $a_{\mathrm{best}}=a_{\kappa_1}$;

  \ELSE  \STATE Increment the counter $I'=I'+1$;
  \ENDIF
  \STATE increment the counter $I=I+1$;
  \ENDWHILE
  \STATE Output the best solution $a_{\mathrm{best}}$ and its value $f(a_{\mathrm{best}})$.
 \end{algorithmic}
  \caption{The proposed generalized extremal optimization procedure.}\label{alg:GEO}
\end{figure}


In the above proposed algorithm, the value of $n_*$ may be chosen differently each time before the set $\Omega_*$ is determined. If $n_*$ is varied cyclicly from 1 through some fixed number, for example, 3, then the resulted algorithm shall be called a modified GEO (MGEO). Specifically, the MGEO firstly generates the neighboring task-to-QC assignments by using $n_*$=1, subsequently carries out the same algorithm but with $n_*$=2, and then with $n_*$=3, and then go back to $n_*$=1 and so on. In the instance that the value of $n_*$ remains unchanged during the computation, the algorithm shall be simply called the GEO.

In contrast to other population-based algorithms,  the proposed GEO/MGEO has only    one string of bits, instead of having a set of strings.   However the proposed algorithm  can search all possible neighboring task-to-QC assignments once the set  $\Omega_*$ of tasks to be mutated are chosen for a given incumbent assignment. In addition, the proposed algorithm updates the incumbent task-to-QC assignment by picking a suitable one from a new set of neighboring assignments of the   incumbent one.  Thus with the continued execution of the searching process, variable neighborhoods in the solution space are searched. Consequently the proposed algorithm may not get trapped in a local optimum, and may find more good solutions,    even more likely to find the global optimal solution.
   In addition, in order to accelerate the convergence of the algorithm, the parameter $\tau$ should not be too small. According to the preliminary tests, when $\tau=5$, the algorithm can obtain  better performance in terms of the solution quality and the computational time. Therefore we choose $\tau=5$ in our experiments when using the  GEO/MGEO.

In order to obtain the best performance of the proposed GEO algorithm,  the issue of selecting a suitable number of simultaneously mutated tasks is studied. Such a suitable number should balance the computing time with optimality of the solution. In the proposed algorithm, the number of simultaneously mutated tasks $n_*$ is set not more than three due to the significantly increased computing time as $n_*$ increases. Consequently, there are three versions of the proposed GEO, each with the value of $n_*$ equal to one of 1, 2 or 3. Meanwhile, a modified GEO (MGEO) algorithm is proposed such that the MGEO performs a sequence of mutation procedures so that $n_*$ takes on cyclic values of 1, 2 and 3. The entire neighboring task-to-QC search process is iterated until the stopping criterion is met.

\begin{figure}[h]
  \centering
  \includegraphics[scale=0.4]{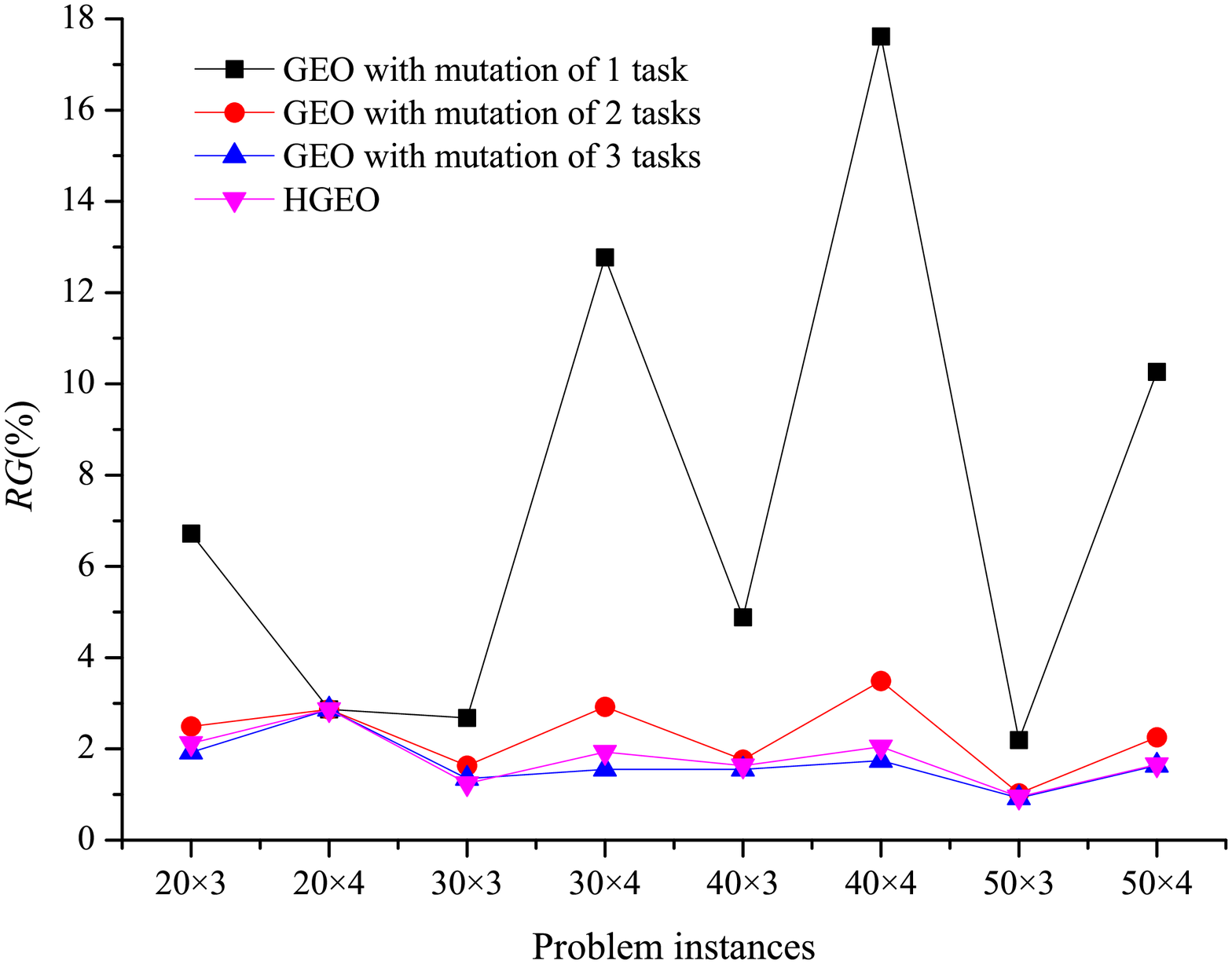}\\
  \caption{Computational results of the four algorithms}\label{fig:RG_C}
\end{figure}

\begin{figure}[h]
  \centering
  \includegraphics[scale=0.4]{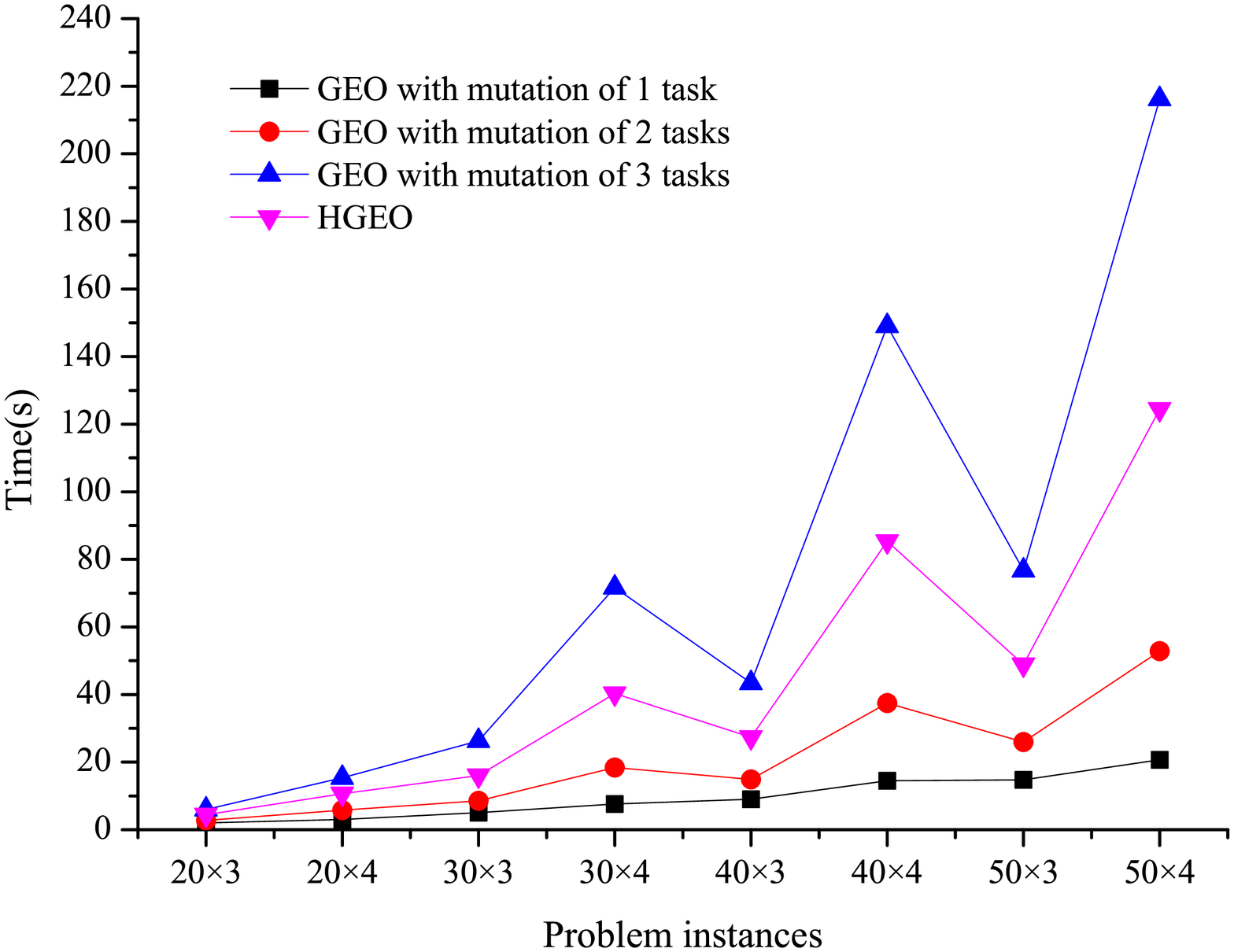}\\
  \caption{Computational times of the four algorithms}\label{fig:Time_C}
\end{figure}

The three versions of the GEO and the MGEO are tested on eight randomly generated instances with $n\in \{20, 30, 40, 50\}$ and $m\in\{3, 4\}$. In these instances, the configurations of producing instances are the same as those of \citet{KimPark2004}. The performance of an algorithm is measured by the relative gap \emph{RG} defined by the equation
\begin{equation}\label{eq:eq_RG}
  RG(\%)=(f_{\mathrm{alg}}-f_{\mathrm{opt}})/f_{\mathrm{opt}} \times 100
\end{equation}
where $f_{\mathrm{alg}}$ is the average solution value of 30 independent runs by a specific  GEO algorithm
 and  $f_{\mathrm{opt}}$  is the optimum makespan. When $f_{\mathrm{opt}}$ is not available for a given instance,  it is replaced by   a lower bound  $f_{\mathrm{LB}}$     of the makespan  that is calculated by the expression
 \begin{equation}\label{eq:eq_LB}
   f_{\mathrm{LB}}\geqslant r_{0}^{k} + \sum_{i\in \Omega}x_{i}^{k}\cdot p_{i} + \max_{i \in \Omega}(t_{0i}^{k}\cdot x_{i}^{k}), \enskip \qquad k\in Q.
 \end{equation}
 In equation \eqref{eq:eq_LB}, the right-hand side gives the shortest service time $c^k$ of a crane $k\in Q$ to complete all assigned tasks without considering any interference constraints. The $f_{\mathrm{LB}}$ is then taken as the maximum of all $c^k$, $k\in Q$. The value of $c^k, \ k\in Q$ is determined by the earliest ready time of crane $k$, the sum of the processing times of tasks assigned to crane $k$ and the maximal traveling time of crane $k$ from its initial bay position to the bay positions of tasks assigned to crane $k$.
 Figure \ref{fig:RG_C} compares the computational results of the four algorithms. Figure \ref{fig:Time_C} shows the run times of those algorithms.
 Note that the values of computing time in figure \ref{fig:Time_C} are the averages of 30 runs for each instance.
 As seen in figure \ref{fig:RG_C}, the MGEO offers about the same quality solutions as those by the GEO with mutation of 3 tasks, and they deliver better solutions (with smaller values of the \emph{RG}s) than those by the other two versions of the GEO. However, the MGEO consumes less time than the GEO with mutation of 3 tasks, as shown in Figure \ref{fig:Time_C}. Therefore, in what follows, the MGEO is adopted to obtain the near-optimal solutions for all test instances.

\section{Numerical experiments}\label{sec:sec5}

In this section, the performances of {\em the proposed approaches (the GEO/MGEO) are evaluated using the benchmark data provided by \citet{KimPark2004}}. Those data have been widely used in the literature. The benchmark data consists of nine instance sets of different problem size with ten different instances each,  as outlined in table \ref{tab:benchmark_Kim}.   In each instance, the processing time for every task is randomly generated from a uniform distribution of $U$(3, 180). It is assumed that the number of ship-bays in the vessel is the same as the number of tasks. The initial locations of cranes are assumed to be equally located on the vessel. The performing order of tasks of a bay is completely determined by the following precedence relation: {discharging tasks on deck, discharging tasks in a hold, loading tasks in a hold, and loading tasks on deck.} {That is to say, when discharging and loading operations must be performed at the same ship-bay, the discharging operation must precede the loading operation;  when a discharging operation is performed in a ship-bay, tasks on a deck must be performed before tasks in the hold of the same ship-bay are performed. Also, the loading operation in a hold must precede the loading operation on the deck of the same ship-bay.}
In addition, the earliest ready time of each crane is set to zero. The traveling time $t_0$ of a crane between adjacent bays is set to one time unit, and the safety margin $\delta$ is set to one bay.

The proposed MGEO algorithm was coded in Matlab 7.11 on a personal computer with a Pentium dual-core 2.6 GHz processor and 2GB RAM. Since the MGEO algorithm is stochastic, the algorithm is run 30 times independently for each instance. In following tables, the listed objective values and the listed times are averages of 30 runs, respectively, unless otherwise indicated.

\begin{table}[h]\scriptsize
  \centering
  \caption{Outlines of the QCSP benchmark data introduced by  \citet{KimPark2004}.}
    \begin{tabular}{rrrr}
    \toprule
    Instance set & Problem no. & Number of tasks & Number of QCs \\
    \midrule
    A     & 13-22 & 10    & 2 \\
    B     & 23-32 & 15    & 2 \\
    C     & 33-42 & 20    & 3 \\
    D     & 43-52 & 25    & 3 \\
    E     & 53-62 & 30    & 4 \\
    F     & 63-72 & 35    & 4 \\
    G     & 73-82 & 40    & 5 \\
    H     & 83-92 & 45    & 5 \\
    I     & 93-102 & 50    & 6 \\
    \bottomrule
    \end{tabular}%
  \label{tab:benchmark_Kim}%
\end{table}%

\begin{figure}[h]
  \centering
  \includegraphics[scale=0.75]{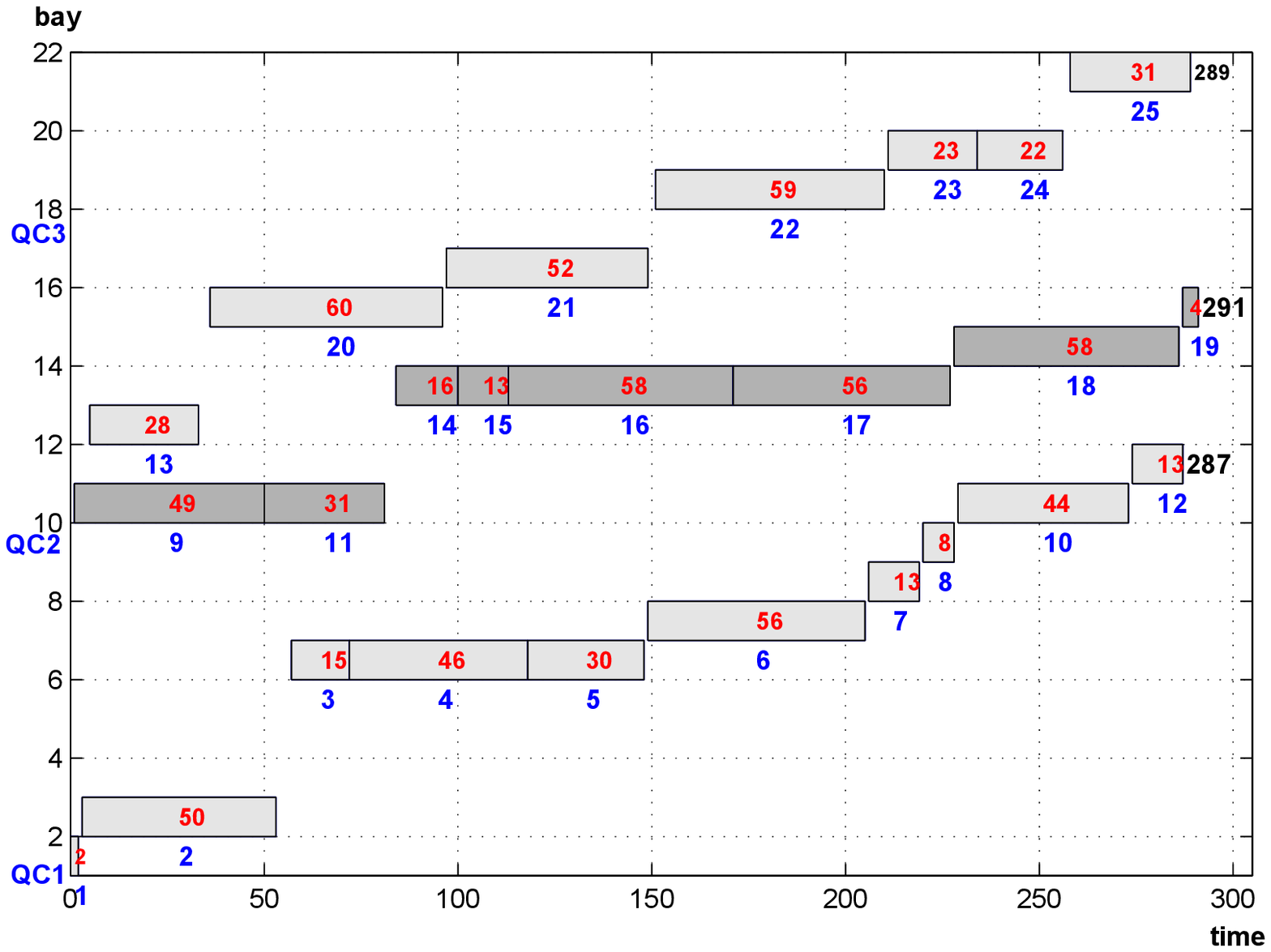}\\
  \caption{The Gantt chart of the schedule  for instance 43 with objective value of 873.}\label{fig:K43_873}
\end{figure}

\begin{figure}[h]
  \centering
  \includegraphics[scale=0.65]{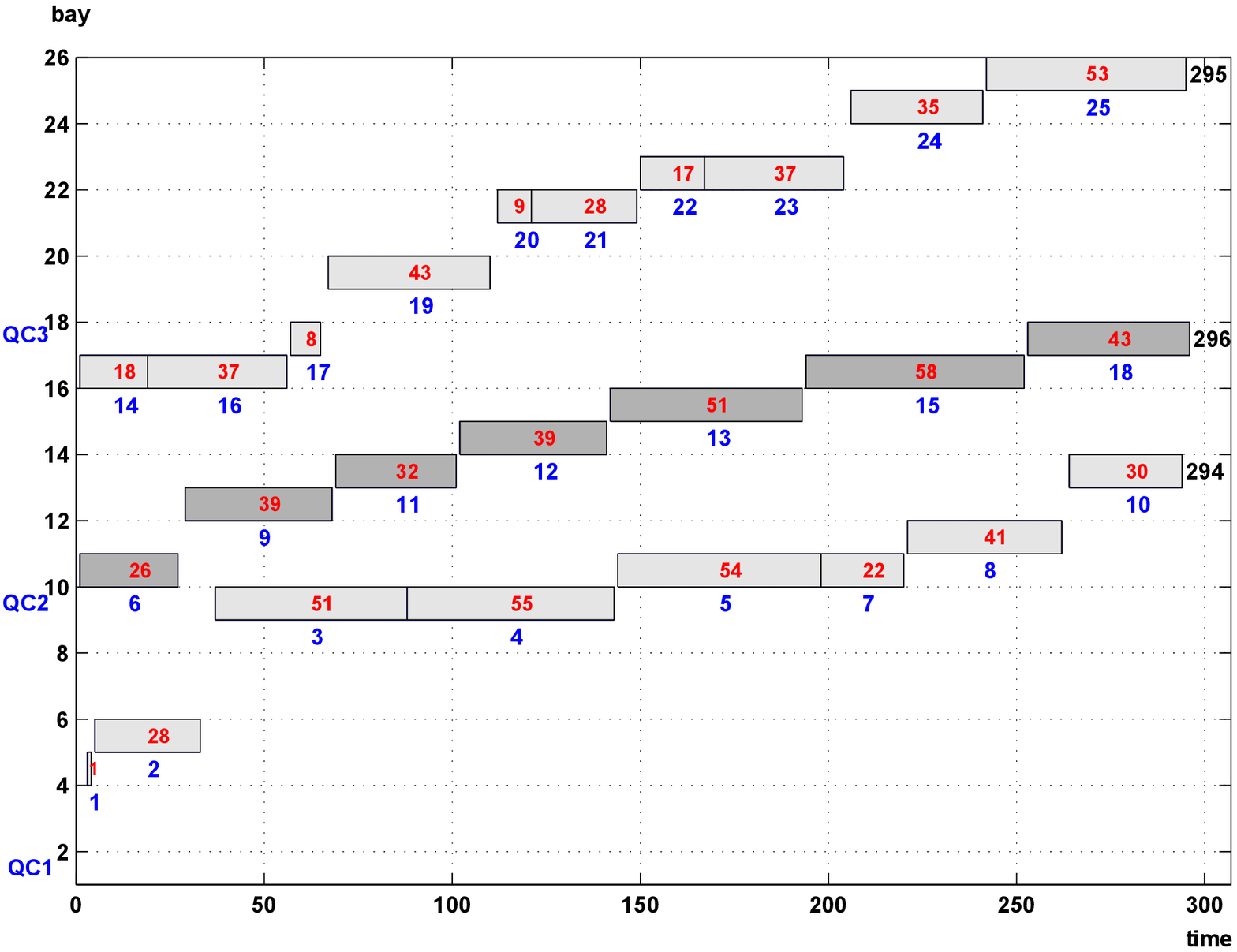}\\
  \caption{The Gantt chart of the schedule  for instance 49 with objective value of 888.}\label{fig:K49_888}
\end{figure}

 {\em The proposed MGEO algorithm was compared with the known methods}, including the Branch-and-Bound algorithm (B\&B) and the GRASP-heuristic of \citet{KimPark2004}, the Branch-and-Cut algorithm (B\&C) of \citet{MocciaCordeau2006}, the Tabu Search (TS) of \citet{SammarraCordeau2007}, the Unidirectional Scheduling heuristic (UDS) of \citet{BierwirthMeisel2009}, and the Modified Genetic Algorithm (MGA) of \citet{Chung20124213}. These approaches have been used to solve a subset of 37 instances from the benchmark problems, labeled 13 to 49.
 The computational results achieved for the 37 instances are listed in table \ref{tab:result_37}.
 {
 If the optimal solution of an instance is not provided in the existing literature or delivered by the CPLEX 12.5 within the 60-minute time limit, a corresponding lower bound is calculated by expression \eqref{eq:eq_LB} for that instance.
 For instances 22 and 42, the two corrected best solutions obtained by \citet{BierwirthMeisel2009} with considering the interference constraints are {\em confirmed} to be {\em optimal} solutions by the calculation of the CPLEX  using the mathematical model presented in Section \ref{sec:Sec3}.
 Moreover,  the best value for instance 39 obtained by known methods is also  {\em confirmed} to be    {\em optimal} by the CPLEX within 47.48 minutes. } For these three instances, their best values were not known to be optimal.  The confirmed optimal values for these three instances are highlighted in the column $f_{\mathrm{opt}}$.
 The results show that {\em the proposed MGEO can obtain the same best known solutions as delivered by other existing optimization approaches in all 37 instances except for instances 19 and 45}. In particular, {\em two new best solutions are found by the MGEO for problem instances 43 and 49}. To verify the two new best solutions, their Gantt charts are shown in figures \ref{fig:K43_873} and \ref{fig:K49_888}. The computational results  indicate that the MGEO is effective in solving  small- and medium-sized instances of the QCSP.
 In addition, {\em two better lower bounds are found by using equation \eqref{eq:eq_LB} compared with that of \citet{BierwirthMeisel2009} in instances 45 and 49}.

\begin{table}[h]\scriptsize
  \centering
  \caption{The results of   computational experiments for  37 instances.}

\begin{tabular}{ccccccccc}
    \toprule
    Problem & \multirow{2}[4]{*}{$f_{\mathrm{opt}}$} & \multicolumn{2}{c}{KP(2004)} & MCGL(2006) & SCLM(2007) & BM(2009) & CC(2012) & \multirow{2}[4]{*}{MGEO} \\
        \cmidrule[0.05em](r){3-4}
    \cmidrule[0.05em](lr){5-5}
    \cmidrule[0.05em](lr){6-6}
    \cmidrule[0.05em](lr){7-7}
    \cmidrule[0.05em](l){8-8}

    no.   &       & B\&B  & GRASP & B\&C  & TS    & UDS   & MGA   &  \\
    \midrule
    13    & 453   & 453   & 453   & 453   & 453   & 453   & 453   & \textbf{453} \\
    14    & 546   & 546   & 546   & 546   & 546   & 546   & 546   & \textbf{546} \\
    15    & 513   & 513   & 516   & 513   & 513   & 513   & 513   & \textbf{513} \\
    16    & 312   & 321   & 321   & 312   & 312   & 312   & 312   & \textbf{312} \\
    17    & 453   & 456   & 456   & 453   & 453   & 453   & 453   & \textbf{453} \\
    18    & 375   & 375   & 375   & 375   & 375   & 375   & 375   & \textbf{375} \\
    19    & 543   & 552   & 552   & 543   & 543   & 543   & 543   & 552 \\
    20    & 399   & 480   & 480   & 399   & 399   & 399   & 399   & \textbf{399} \\
    21    & 465   & 465   & 465   & 465   & 465   & 465   & 465   & \textbf{465} \\
    22    & \textbf{540}   & 720   & 720   & 537   & 537   & 540   & 537   & \textbf{540} \\
    23    & 576   & 576   & 591   & 576   & 582   & 576   & 576   & \textbf{576} \\
    24    & 666   & 669   & 675   & 666   & 669   & 666   & 669   & \textbf{666} \\
    25    & 738   & 738   & 741   & 738   & 741   & 738   & 744   & \textbf{738} \\
    26    & 639   & 639   & 651   & 639   & 639   & 639   & 645   & \textbf{639} \\
    27    & 657   & 657   & 687   & 660   & 660   & 657   & 660   & \textbf{657} \\
    28    & 531   & 537   & 549   & 531   & 531   & 531   & 531   & \textbf{531} \\
    29    & 807   & 807   & 819   & 807   & 810   & 807   & 810   & \textbf{807} \\
    30    & 891   & 891   & 906   & 891   & 891   & 891   & 897   & \textbf{891} \\
    31    & 570   & 570   & 570   & 570   & 570   & 570   & 570   & \textbf{570} \\
    32    & 591   & 591   & 597   & 591   & 591   & 591   & 594   & \textbf{591} \\
    33    & 603   & 603   & 666   & 603   & 603   & 603   & 603   & \textbf{603} \\
    34    & 717   & 717   & 762   & 717   & 735   & 717   & 717   & \textbf{717} \\
    35    & 684   & 678-690 & 699   & 684   & 690   & 684   & 690   & \textbf{684} \\
    36    & 678   & 678-720 & 708   & 678   & 681   & 678   & 636   & \textbf{678} \\
    37    & 510   & 510-516 & 540   & 510   & 519   & 510   & 522   & \textbf{510} \\
    38    & $\text{606}^*$   & 609-633 & 660   & 618   & 618   & 618   & 618   & \textbf{618} \\
    39    & \textbf{513}   & 513-552 & 579   & 513   & 519   & 513   & 519   & \textbf{513} \\
    40    & 564   & 558-576 & 597   & 564   & 567   & 564   & 567   & \textbf{564} \\
    41    & $\text{582}^*$   & 585-654 & 642   & 588   & 594   & 588   & 588   & \textbf{588} \\
    42    & \textbf{573}   & 561-588 & 666   & 570   & 576   & 573   & 576   & \textbf{573} \\
    43    & $\text{858}^*$ & 864-951 & 942   & 897   & 879   & 876   & 897   & \text{\textbf{873}}$^\diamond$ \\
    44    & $\text{813}^*$ & 813-879 & 858   & 822   & 834   & 822   & 855   & \textbf{822} \\
    45    & $\text{825}^*$   & 825-861 & 873   & 840   & 852   & 834   & 864   & 837 \\
    46    & 690   & 687-708 & 735   & 690   & 690   & 690   & 723   & \textbf{690} \\
    47    & 792   & 789-912 & 807   & 792   & 792   & 792   & 819   & \textbf{792} \\
    48    & $\text{618}^*$ & 627-669 & 669   & 645   & 663   & 639   & 663   & \textbf{639} \\
    49    & $\text{885}^*$   & 888-915 & 972   & 927   & 912   & 894   & 915   & \text{\textbf{888}}$^\diamond$ \\
    \bottomrule

    \multicolumn{3}{l}{$^*$  lower bound}\\
    \multicolumn{3}{l}{$^\diamond$ new best solution}

    \end{tabular}%
  \label{tab:result_37}%
\end{table}%

The average computational time measured in minutes for each of the known methods is reported in table \ref{tab:time_1349}, together with the specification of the PC used.
It is found that the MGEO   performs very well compared with other approaches. Note that the UDS heuristic consumes slightly  less time than the MGEO.

\begin{table}[h]\scriptsize
  \centering
  \caption{Computational times of existing approaches and the MGEO.}
    \begin{tabular}{cccccccc}
    \toprule
    \multirow{3}[4]{*}{Problem set} & \multicolumn{2}{c}{KP(2004)} & MCGL(2006) & SCLM(2007) & BM(2009) & CC(2012) &  \\

            & \multicolumn{2}{c}{P2(466MHZ)\_64MB} & P4(2.5GHz)\_512MB  & P4(2.5GHz)\_512MB    & P4(2.8GHz)   & Core2(2.0GHz)\_2GB   & Core2(2.6GHz)\_2GB\\

           \cmidrule[0.05em](r){2-3}
    \cmidrule[0.05em](lr){4-4}
    \cmidrule[0.05em](lr){5-5}
    \cmidrule[0.05em](l){6-6}
    \cmidrule[0.05em](l){7-7}
    \cmidrule[0.05em](l){8-8}
          & B\&B  & GRASP & B\&C  & TS    & UDS   & MGA   & MGEO \\
           \midrule
    A     & 0.44  & 0.35  & 1.01  & 1.52  & 1.12$\times$10$^{-5}$ & 0.52  & 8.08$\times$10$^{-4}$ \\
    B     & 17.53  & 1.46  & 8.91  & 5.86  & 3.68$\times$10$^{-5}$ & 0.75  & 3.11$\times$10$^{-3}$ \\
    C     & 564.47  & 3.16  & 72.19  & 21.75  & 6.26$\times$10$^{-4}$ & 1.18  & 3.34$\times$10$^{-2}$ \\
    D     & 750.87  & 7.56  & 102.49  & 48.68  & 3.43$\times$10$^{-3}$ & 1.58  & 6.63$\times$10$^{-2}$ \\
    \bottomrule

    \end{tabular}%
  \label{tab:time_1349}%
\end{table}%

\begin{figure}[h]
  \centering
  \includegraphics[scale=0.4]{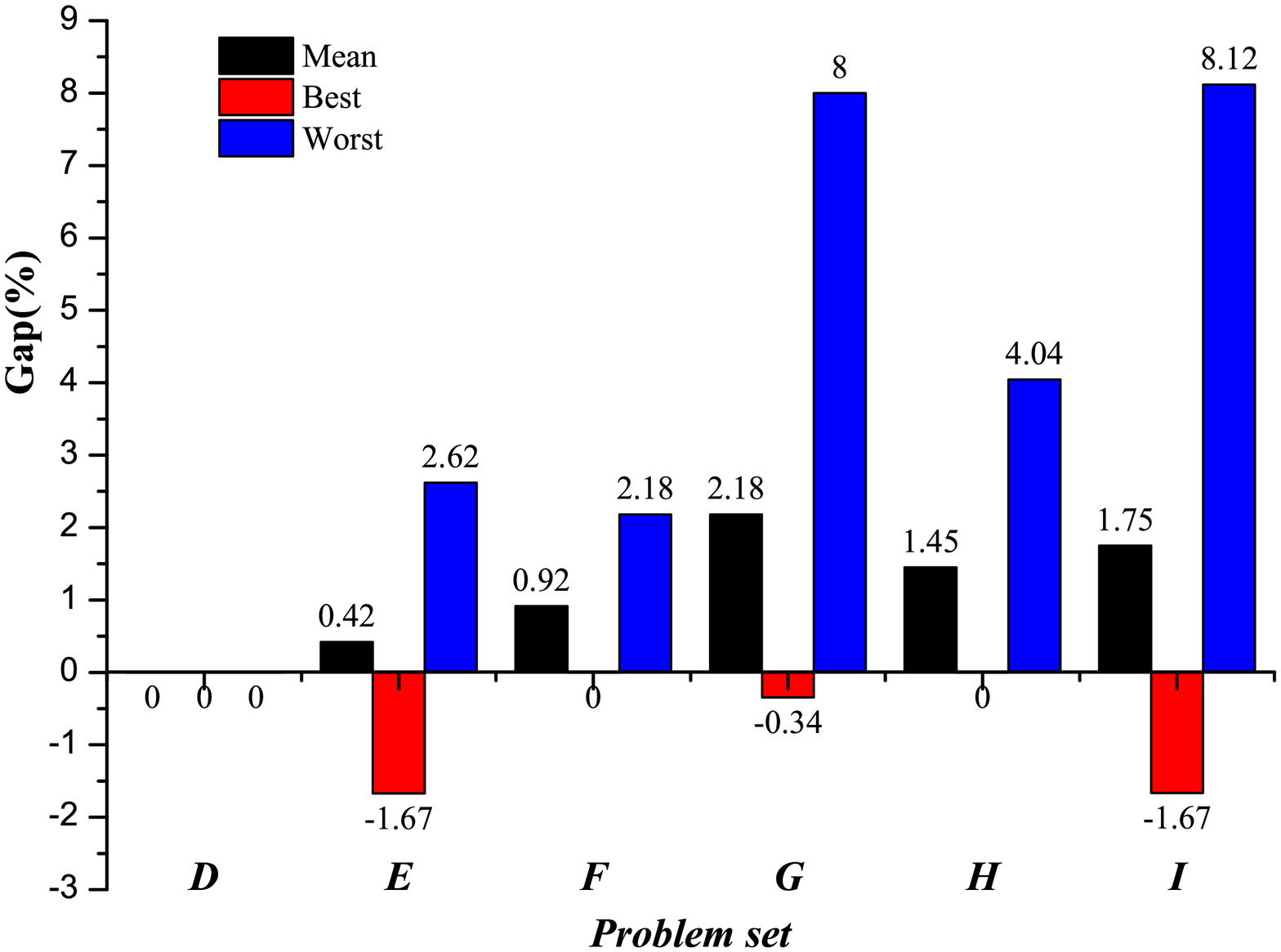}\\
  \caption{Comparison of the average solution performance between the MGEO and the UDS.}\label{fig:DIGAP}
\end{figure}

Based on the good performance of the MGEO on  instances 13-49, the proposed algorithm is used to handle  larger instances provided in sets E to I of the benchmark problems. These instances are only solved by \citet{BierwirthMeisel2009}, the MGEO is compared with the UDS heuristic. Since the MGEO algorithm is stochastic, for each instance, the mean solution value, the best solution value and the worst solution of 30 runs are listed  in table \ref{tab:50102}. Among these instances, the largest size has 50 tasks and 6 QCs.  It is found that the solution quality achieved by the MGEO is very good. The average\emph{ RG}$_\mathrm{Mean}$ for the 53 instances is 3.34\% as   shown in the last row of table \ref{tab:50102}, which is only slightly higher than the average \emph{RG} 2.04\% given by the UDS. The average  RG of the best solutions given by the MGEO is 2.43\%, which is very close to the one of the UDS.
Nevertheless, the MGEO algorithm delivers solutions for these large-sized instances very quickly. The average time of the MGEO is merely 1.36 minutes. Even for the most intractable instance 97, the MGEO only consumes 5.36 minutes to solve it. In addition, {\em six new best solutions are found by the MGEO} as highlighted in the column``Best" under ``MGEO" in table \ref{tab:50102}. When the  size of instances increases, the computational time of the  MGEO only increase slightly.
In order to evaluate the performance of the MGEO and the UDS, the average-in-set gaps are showed in figure \ref{fig:DIGAP}. In figure \ref{fig:DIGAP}, the gap is calculated in percent of the solution $f_{\mathrm{MGEO}}$ against $f_{\mathrm{UDS}}$, i.e. gap=$(f_{\mathrm{MGEO}}-f_{\mathrm{UDS}})/f_{\mathrm{UDS}}\times 100$.

Summarizing the findings of the above computational tests, the developed MGEO algorithm is capable of solving all instances in a reasonable time in comparison with the existing QCSP solution methods.
It delivers high quality solutions for small to medium size instances and performs as well  for large-sized problems.
The computational times of the MGEO are short, especially for large-sized instances. Therefore, the MGEO is robust with respect to both heuristics and meta-heuristics in terms of solution quality and computational time.


\begin{table}[h]\scriptsize
  \centering
  \caption{Comparison between MGEO and UDS for large-sized instances.}
    \begin{tabular}{cccccccccccc}
    \toprule
    Problem & \multirow{2}[0]{*}{LB} & \multicolumn{3}{c}{UDS} & \multicolumn{7}{c}{MGEO}                      \\
          \cmidrule[0.05em](r){3-5}
    \cmidrule[0.05em](l){6-12}

    no.   &       & Best  & $RG_{\mathrm{Best}}$ & Time(m) & Mean  & Best  & Worst & $RG_{\mathrm{Mean}}$ & $RG_{\mathrm{Best}}$ & $RG_{\mathrm{Worst}}$ & Time(m) \\
     \midrule
    50    & 738   & 741   & 0.41  & 0.03  & 741   & 741   & 741   & 0.41  & 0.41  & 0.41  & 0.05  \\
    51    & 780   & 798   & 2.31  & 0.01  & 798   & 798   & 798   & 2.31  & 2.31  & 2.31  & 0.06  \\
    52    & 948   & 960   & 1.27  & 0.07  & 960   & 960   & 960   & 1.27  & 1.27  & 1.27  & 0.10  \\
    53    & 672   & 717   & 6.70  & 60.00  & 710.4 & \textbf{705}$^\diamond$   & 717   & 5.71  & 4.91  & 6.70  & 0.26  \\
    54    & 762   & 774   & 1.57  & 0.02  & 778.1 & 774   & 780   & 2.11  & 1.57  & 2.36  & 0.17  \\
    55    & 672   & 684   & 1.79  & 0.01  & 685.6 & 684   & 693   & 2.02  & 1.79  & 3.13  & 0.20  \\
    56    & 681   & 690   & 1.32  & 0.22  & 692.8 & 690   & 699   & 1.73  & 1.32  & 2.64  & 0.22  \\
    57    & 687   & 705   & 2.62  & 0.24  & 707.9 & 705   & 711   & 3.04  & 2.62  & 3.49  & 0.19  \\
    58    & 777   & 786   & 1.16  & 0.17  & 788   & \textbf{783}$^\diamond$   & 798   & 1.42  & 0.77  & 2.70  & 0.22  \\
    59    & 678   & 687   & 1.33  & 0.01  & 691   & 687   & 705   & 1.92  & 1.33  & 3.98  & 0.32  \\
    60    & 771   & 783   & 1.56  & 0.19  & 790.9 & 786   & 798   & 2.58  & 1.95  & 3.50  & 0.27  \\
    61    & 621   & 639   & 2.90  & 0.04  & 644   & 639   & 654   & 3.70  & 2.90  & 5.31  & 0.17  \\
    62    & 837   & 837   & 0.00  & 0.01  & 844.65 & 837   & 858   & 0.91  & 0.00  & 2.51  & 0.26  \\
    63    & 939   & 948   & 0.96  & 1.51  & 954   & 948   & 957   & 1.60  & 0.96  & 1.92  & 0.38  \\
    64    & 729   & 741   & 1.65  & 1.06  & 744.9 & 741   & 750   & 2.18  & 1.65  & 2.88  & 0.43  \\
    65    & 837   & 837   & 0.00  & 1.61  & 844.5 & 840   & 849   & 0.90  & 0.36  & 1.43  & 0.37  \\
    66    & 921   & 924   & 0.33  & 0.63  & 933.6 & 930   & 939   & 1.37  & 0.98  & 1.95  & 0.49  \\
    67    & 864   & 882   & 2.08  & 0.24  & 889.8 & 885   & 897   & 2.99  & 2.43  & 3.82  & 0.37  \\
    68    & 960   & 963   & 0.31  & 0.03  & 975   & 963   & 984   & 1.56  & 0.31  & 2.50  & 0.29  \\
    69    & 792   & 807   & 1.89  & 1.40  & 811.5 & 810   & 816   & 2.46  & 2.27  & 3.03  & 0.32  \\
    70    & 945   & 957   & 1.27  & 0.61  & 965.4 & 957   & 966   & 2.16  & 1.27  & 2.22  & 0.26  \\
    71    & 831   & 834   & 0.36  & 3.77  & 845.4 & 837   & 849   & 1.73  & 0.72  & 2.17  & 0.39  \\
    72    & 723   & 744   & 2.90  & 0.35  & 752.4 & 744   & 756   & 4.07  & 2.90  & 4.56  & 0.32  \\
    73    & 849   & 870   & 2.47  & 31.71  & 882.9 & \textbf{867}$^\diamond$   & 888   & 3.99  & 2.12  & 4.59  & 1.42  \\
    74    & 834   & 843   & 1.08  & 4.71  & 854.1 & 849   & 864   & 2.41  & 1.80  & 3.60  & 0.62  \\
    75    & 672   & 675   & 0.45  & 0.37  & 704.4 & 687   & 729   & 4.82  & 2.23  & 8.48  & 1.17  \\
    76    & 828   & 852   & 2.90  & 0.90  & 870.9 & 861   & 885   & 5.18  & 3.99  & 6.88  & 0.82  \\
    77    & 678   & 699   & 3.10  & 1.27  & 717.6 & 708   & 729   & 5.84  & 4.42  & 7.52  & 1.34  \\
    78    & 630   & 642   & 1.90  & 8.96  & 648.6 & 645   & 651   & 2.95  & 2.38  & 3.33  & 0.82  \\
    79    & 726   & 744   & 2.48  & 1.52  & 772.8 & 768   & 780   & 6.45  & 5.79  & 7.44  & 1.20  \\
    80    & 726   & 750   & 3.31  & 1.28  & 764.4 & 756   & 774   & 5.29  & 4.13  & 6.61  & 0.87  \\
    81    & 711   & 738   & 3.80  & 1.28  & 747.6 & 744   & 762   & 5.15  & 4.64  & 7.17  & 1.00  \\
    82    & 705   & 717   & 1.70  & 1.03  & 729   & 717   & 735   & 3.40  & 1.70  & 4.26  & 1.40  \\
    83    & 930   & 948   & 1.94  & 6.37  & 973.2 & 954   & 984   & 4.65  & 2.58  & 5.81  & 1.56  \\
    84    & 888   & 897   & 1.01  & 3.29  & 907.8 & 897   & 927   & 2.23  & 1.01  & 4.39  & 2.34  \\
    85    & 954   & 972   & 1.89  & 5.82  & 979.2 & 972   & 984   & 2.64  & 1.89  & 3.14  & 1.40  \\
    86    & 795   & 816   & 2.64  & 60.00  & 829.8 & 822   & 849   & 4.38  & 3.40  & 6.79  & 1.89  \\
    87    & 855   & 867   & 1.40  & 60.00  & 877.2 & 870   & 885   & 2.60  & 1.75  & 3.51  & 2.44  \\
    88    & 747   & 768   & 2.81  & 43.73  & 778.2 & 771   & 786   & 4.18  & 3.21  & 5.22  & 1.18  \\
    89    & 831   & 843   & 1.44  & 10.96  & 855   & 843   & 861   & 2.89  & 1.44  & 3.61  & 1.12  \\
    90    & 1032  & 1053  & 2.03  & 24.95  & 1072.8 & 1053  & 1086  & 3.95  & 2.03  & 5.23  & 2.39  \\
    91    & 822   & 837   & 1.82  & 10.74  & 837   & 837   & 837   & 1.82  & 1.82  & 1.82  & 1.41  \\
    92    & 882   & 897   & 1.70  & 34.61  & 918.6 & 912   & 927   & 4.15  & 3.40  & 5.10  & 1.43  \\
    93    & 795   & 816   & 2.64  & 60.00  & 824.4 & 816   & 834   & 3.70  & 2.64  & 4.91  & 3.78  \\
    94    & 780   & 786   & 0.77  & 60.00  & 805.2 & 798   & 810   & 3.23  & 2.31  & 3.85  & 3.45  \\
    95    & 804   & 834   & 3.73  & 60.00  & 847.8 & 840   & 852   & 5.45  & 4.48  & 5.97  & 4.49  \\
    96    & 795   & 819   & 3.02  & 60.00  & 822.6 & \textbf{816}$^\diamond$   & 828   & 3.47  & 2.64  & 4.15  & 4.02  \\
    97    & 702   & 720   & 2.56  & 60.00  & 726   & \textbf{717}$^\diamond$   & 729   & 3.42  & 2.14  & 3.85  & 5.36  \\
    98    & 720   & 735   & 2.08  & 23.79  & 754.2 & 747   & 765   & 4.75  & 3.75  & 6.25  & 3.37  \\
    99    & 825   & 852   & 3.27  & 60.00  & 858   & 855   & 861   & 4.00  & 3.64  & 4.36  & 4.37  \\
    100   & 867   & 900   & 3.81  & 60.00  & 894.6 & \textbf{885 }$^\diamond$  & 909   & 3.18  & 2.08  & 4.84  & 3.04  \\
    101   & 774   & 813   & 5.04  & 60.00  & 856.2 & 828   & 879   & 10.62  & 6.98  & 13.57  & 2.56  \\
    102   & 879   & 903   & 2.73  & 60.00  & 930.6 & 927   & 933   & 5.87  & 5.46  & 6.14  & 3.83  \\
    \multicolumn{2}{c}{average }       &       & 2.04  & 17.92  &       &       &       & 3.34  & 2.43  & 4.32  & 1.36  \\
    \bottomrule
    \\
    \multicolumn{3}{l}{$^\diamond$ new best solution}
    \end{tabular}%
  \label{tab:50102}%
\end{table}%

\section{Conclusions}\label{sec:sec6}

In this article, the QCSP at a port container terminal is studied while interference constraints such as no-crossing and safety margin constraints are considered. A stochastic algorithm called the MGEO is designed to solve the problem. Computational results for a set of well known benchmark problems are performed to evaluate the performance of the proposed algorithm.
By comparing the results delivered by the proposed MGEO and the ones obtained by other existing approaches including B\&B, GRASP, B\&C, TS, UDS, and MGA, it is found that the proposed MGEO performs robust in small-, medium- and large-sized instances and obtains new better solutions in some instances. Besides, the computational time required by the MGEO is much less than the current known approaches in small- and medium- sized instances except for the UDS. Moreover, the computational  time of the MGEO is clearly shorter than that of the UDS in large-sized instances.

For future work, the QC-to-vessel assignment problem can be integrated into the QCSP. Such a more comprehensive model considering both QC-to-vessel assignment and task-to-QC scheduling is expected to provide a better operation plan for a container terminal. On the other hand, the conditions under which a unidirectional search leading to an optimal task-to-QC schedule may be of interest.

\section*{Acknowledgments}

This work was partially supported by the National Natural Science Foundation of China (No. 51175442), the Youth Foundation for Humanities and Social Sciences of Ministry of Education of China (No. 12YJCZH296) and the Fundamental Research Funds for the Central Universities (No.SWJTU09CX022; 2010ZT03).
{  The authors would like to thank Drs. K.H. Kim from Pusan National  University and Y.M. Park from Korea Naval  Academy, who kindly provided the benchmark data suites. We are grateful to Dr. Frank Meisel from Martin-Luther University for his invaluable assistance to our questions. Thanks are also due to referees for the constructive comments on improving our article.  }


\bibliographystyle{gENO}
\addcontentsline{toc}{section}{\refname}
\bibliography{GEO_QCSP72}

\label{lastpage}

\end{document}